\documentclass[10pt]{article}
\usepackage{pb-diagram}
\usepackage{amssymb,latexsym}
\usepackage{amsfonts}
\usepackage{amsmath,amsthm,graphicx}

\newcommand{\bd}{\begin{document}}
\newcommand{\ed}{\end{document}}
\newcommand{\bc}{\begin{center}}
\newcommand{\ec}{\end{center}}
\newcommand{\vs}{\vspace}
\newcommand{\hs}{\hspace}
\newcommand{\bq}{\begin{quote}}
\newcommand{\eq}{\end{quote}}
\newcommand{\mb}{\makebox}
\newcommand{\lt}{\left}
\newcommand{\rt}{\right}
\newcommand{\beqa}{\begin{eqnarray*}}
\newcommand{\eeqa}{\end{eqnarray*}}
\newcommand{\beqn}{\begin{eqnarray}}
\newcommand{\eeqn}{\end{eqnarray}}
\newcommand{\bbibl}{\begin{thebibliography}{99}}
\newcommand{\ebibl}{\end{thebibliography}}
\newcommand{\ti}{\times}
\newcommand{\bit}{\begin{itemize}}
\newcommand{\eit}{\end{itemize}}
\newcommand{\ben}{\begin{enumerate}}
\newcommand{\een}{\end{enumerate}}
\newcommand{\lb}{\label}

\newcommand{\hf}{\hspace*{\fill}}
\newcommand{\vf}{\vspace*{\fill}}
\newcommand{\beq}{\begin{equation}}
\newcommand{\eeq}{\end{equation}}
\newcommand{\ba}{\begin{array}}
\newcommand{\ea}{\end{array}}
\newcommand{\del}{\partial}
\newcommand{\ot}{\otimes}
\newcommand{\nn}{\nonumber}
\newcommand{\R}{\mb{$I\!\!R$}}
\newcommand{\C}{{\cal C}}
\newcommand{\M}{{\cal M}}
\newcommand{\E}{{\cal E}}
\newcommand{\N}{{\cal N}}
\newcommand{\B}{{\cal B}}
\newcommand{\X}{{\cal X}}
\newcommand{\Y}{{\cal Y}}
\newcommand{\F}{{\cal F}}
\newcommand{\Rc}{{\cal R}}
\newcommand{\A}{{\cal A}}
\renewcommand{\P}{{\cal P}}
\renewcommand{\S}{{\cal S}}
\newcommand{\es}{\emptyset}
\newcommand{\ci}{\subseteq}
\newcommand{\cs}{\supseteq}
\renewcommand{\u}{\cup}
\renewcommand{\i}{\cap}
\newcommand{\bu}{\bigcup}
\newcommand{\bi}{\bigcap}
\newcommand{\la}{\leftarrow}
\newcommand{\ra}{\rightarrow}
\newcommand{\Ra}{\Rightarrow}
\newcommand{\Lra}{\Leftrightarrow}
\newcommand{\lgra}{\longrightarrow}
\newcommand{\Lgra}{\Longrightarrow}
\newcommand{\lglra}{\longleftrightarrow}
\newcommand{\Lglra}{\Longleftrightarrow}
\renewcommand{\a}{\alpha}
\renewcommand{\b}{\beta}
\newcommand{\g}{\gamma}
\newcommand{\G}{\Gamma}
\renewcommand{\d}{\delta}
\newcommand{\D}{\Delta}
\newcommand{\e}{\varepsilon}
\newcommand{\eps}{\epsilon}
\newcommand{\h}{\eta}
\renewcommand{\l}{\lambda}
\newcommand{\m}{\mu}
\newcommand{\n}{\nu}
\newcommand{\p}{\pi}
\newcommand{\s}{\sigma}
\newcommand{\Si}{\Sigma}
\newcommand{\ta}{\tau}
\newcommand{\ph}{\phi}
\newcommand{\Ph}{\Phi}
\renewcommand{\c}{\chi}
\newcommand{\om}{\omega}
\newcommand{\Om}{\Omega}
\newcommand{\tri}{\triangle}
\newcommand{\rec}[1]{\frac{1}{#1}}
\newcommand{\sm}[2]{\sum_{#1}^{#2}}
\newcommand{\ld}{\ldots}
\newcommand{\ov}{\overline}
\newcommand{\ol}[1]{$\bar{\mb{#1}}$}
\newcommand{\un}{\underline}
\newcommand{\iy}{\infty}

\newcommand{\wt}{\widetilde}
\newcommand{\ds}{\displaystyle}
\newcommand{\bdm}{\begin{displaymath}}
\newcommand{\edm}{\end{displaymath}}

\newcommand{\nin}{\not\in}
\newcommand{\bt}{\begin{tabular}}
\newcommand{\et}{\end{tabular}}
\newcommand{\alter}[2]{\lt\{ \ba {ll}#1 \\ #2 \ea \rt.}
\newcommand{\alt}[4]{\lt\{ \ba{ll}#1 & \mb{if \,\,}#2 \\ #3 & \mb{if
               \,\,}#4 \ea \rt.}
\newcommand{\alto}[6]{ \lt\{ \ba{ll}#1 & \mb{if \,\,}#2 \\ #3 & \mb{if
               \,\,} #4 \\ #5 & \mb{if \,\,}#6 \ea \rt.}
\newcommand{\altero}[5]{\mb{$\lt\{ \ba {ll}#1 & \mb{if \,\,}#2 \\ #3 &
               \mb{if \,\,} #4 \\ #5 & \mb{otherwise} \ea \rt.$}}

\newcounter{cnt1}
\newcounter{cnt2}
\newcounter{cnt3}
\newcommand{\blr}{\begin{list}{$($\roman{cnt1}$)$} {\usecounter{cnt1}
        \setlength{\topsep}{0pt} \setlength{\itemsep}{0pt}}}
\newcommand{\bla}{\begin{list}{$($\alph{cnt2}$)$} {\usecounter{cnt2}
        \setlength{\topsep}{0pt} \setlength{\itemsep}{0pt}}}
\newcommand{\bln}{\begin{list}{$($\arabic{cnt3}$)$} {\usecounter{cnt3}
                \setlength{\topsep}{0pt} \setlength{\itemsep}{0pt}}}
\newcommand{\el}{\end{list}}
\newcommand{\no}{\noindent}
\newtheorem{Thm}{Theorem}[section]
\newtheorem{Lem}[Thm]{Lemma}
\newtheorem{Prop}[Thm]{Proposition}
\newtheorem{Def}[Thm]{Definition}
\newtheorem{Exm}[Thm]{Example}
\newtheorem{Rem}[Thm]{Remark}
\newtheorem{Cor}[Thm]{Corollory}
\renewcommand{\baselinestretch}{1}
\newcommand{\ilim}{\mathop{\varprojlim}\limits}
\newcommand{\dlim}{\mathop{\varinjlim}\limits}

\begin{document}
\title{Steenrod's reduced power operations in simplicial Bredon-Illman cohomology with local coefficients}
\author{Goutam Mukherjee and Debasis Sen}
\date{}
\maketitle
\noindent
\begin{abstract}
In this paper we use Peter May's algebraic approach to Steenrod operations to construct Steenrod's reduced power operations in simplicial Bredon-Illman cohomology with local coefficients of a one vertex $G$-Kan complex, $G$ being a discrete group.
\end{abstract}
{\bf Keywords Simplicial sets, group action, equivariant cohomology, equivariant local coefficients, cohomology operations, Steenrod's reduced power operations.}
\footnote {\bf The second author would like to thank CSIR for its support.}

\section{Introduction}

The study of cohomology operations has been one of the important areas of research in algebraic topology for a long time. For instance, they have been extensively used
to compute obstructions \cite{steen3}, to study homotopy type of complexes \cite{thom} and to show essentiality of maps of spheres \cite{bs}. A class of basic operations are Steenrod's squares and reduced power operations \cite{steen4}, \cite{steen5}, \cite{ara}. Steenrod's squares are defined for cohomology with $\mathbb{Z}_2$ coefficients whereas Steenrod's reduced powers are defined in cohomology with coefficients in $\mathbb{Z}_p$, $p\neq 2$ a prime. A very general and useful method of constructing these operations is given in \cite{may1}. A categorical approach to Steenrod operations can be found in \cite{eps}. In \cite{git}, S. Gitler constructed reduced power operations in cohomology with local coefficients. A well known result of Eilenberg describes cohomology of a space with local coefficients by the cohomology of an invariant subcomplex of its universal cover, equipped with the action of the fundamental group of the space \cite{spa}. The main idea of Gitler's construction is to lift power operations in this invariant cochain subcomplex and reproduce the operations in cohomology with local coefficients via Eilenberg's description. The relevant local coefficients in this context is obtained by a fixed action of the fundamental group of the space on a fixed cyclic group of prime order $p \neq 2.$

Recently, in \cite{ms}, we introduced  simplicial equivariant cohomology with local coefficients, which is the simplicial version of Bredon-Illman cohomology with local coefficients \cite{mm}. The aim of this paper is to construct Steenrod's reduced power operations in simplicial Bredon-Illman cohomology with local coefficients, where the equivariant local coefficients take values in a $\mathbb{Z}_p$-algebra, for a prime $p>2$. Throughout our method is simplicial. It may be mentioned that in \cite{ginot}, for a space with a group action Steenrod's squares have been introduced in Bredon-Illman cohomology with local coefficients.

We have the notion of `universal $O_G$-covering complex' of a one vertex $G$-Kan complex $X$ \cite{ms}. This is defined as a contravariant functor from the category of canonical orbits to the category of one vertex Kan complexes and is the analogue, in the equivariant context, of the universal cover of a one vertex Kan complex \cite{gugg}. This universal $O_G$-covering complex comes equipped with an action of an $O_G$-group $\underline{\pi}X$ (see Section $4$ for details) and an equivariant analogue of Eilenberg theorem holds \cite{ms1}. Following Gitler \cite{git}, we first construct the power operations in the $\underline{\pi}X$-equivariant cohomology of `universal $O_G$-covering complex'. This is done by applying the algebraic description of Steenrod power operations of P. May \cite{may1}. We then use the equivariant version of Eilenberg theorem to reproduce Steenrod's reduced power operations in the present context. It may be remarked that our method also applies when $p=2$, and hence yields Steenrod squares too.

The paper is organized as follows. In Section $2$, we recall some standard results and fix notations. The notion of equivariant local coefficients of a simplicial set equipped with a simplicial group action is based on fundamental groupoid. In Section $3$, we recall these concepts and quickly review the definition of simplicial Bredon-Illman cohomology with local coefficients. In Section $4$, we state the equivariant version of Eilenberg theorem. In Section $5$, we briefly recall the algebraic method of P. May and then apply it to construct Steenrod's reduced power operations in simplicial Bredon-Illman cohomology with local coefficients.

\section{Preliminaries}

In this section we set up our notations and recall some standard facts \cite{may}, \cite{gj}.

Throughout $\mathcal{S}$ will denote the category of simplicial sets and simplicial maps.  Let $\Delta[n]$ denote the standard simplicial $n$-simplex and $\Delta_n$ be the unique non-degenerate $n$-simplex of $\Delta[n].$
We have simplicial maps $\delta_i\colon\Delta[n-1]\rightarrow \Delta[n]$ and $\sigma_i\colon\Delta[n+1]\rightarrow \Delta[n]$ for $0\leq i \leq n$ defined by $\delta_i(\Delta_{n-1})=\partial_i(\Delta_n)$ and $\sigma_i(\Delta_{n+1})=s_i(\Delta_n)$. The boundary subcomplex $\partial \Delta[n]$ of $\Delta[n]$ is defined as the smallest subcomplex of $\Delta[n]$ containing the faces $\partial_i\Delta_n,~~ i=0,1,\cdots,n $.

\begin{Def}
Let $G$ be a discrete group.  A $G$-simplicial set is a simplicial object in the category of $G$-sets. More precisely, a $G$-simplicial set is a simplicial set $\{X_n ; \partial_i, s_i, 0\leq i \leq n\}_{n\geq 0}$ such that each $X_n$ is a $G$-set and the face maps $\partial_i\colon X_n\rightarrow X_{n-1}$ and the degeneracy maps $s_i\colon X_n\rightarrow X_{n+1}$ commute with the $G$-action.
A $G$-simplicial map between G-simplicial sets is a simplicial map which commutes with the G-action.
\end{Def}

For a $G$ simplicial set $X$, we consider $X\times \Delta[1]$ as a $G$-simplicial set with trivial $G$-action on $\Delta[1]$.
\begin{Def}
Two $G$-simplicial maps $f,g\colon X\rightarrow Y$ between $G$-simplicial sets $X$ and $Y$ are $G$-homotopic if there exists a $G$-simplicial map $\mathcal{H}\colon X\times \Delta[1]\rightarrow Y$ such that $$\mathcal{H}\circ (id\times \delta_1) =f,~~\mathcal{H}\circ (id\times \delta_0) =g,$$ where $X\times \Delta[0]$ is identified with $X$.
The map $\mathcal{H}$ is called a $G$-homotopy from $f$ to $g$ and we write $\mathcal{H}\colon f\simeq_G g.$
If $i\colon X^{\prime}\subseteq X$ is an inclusion of a subcomplex and $f,~g$ agree on $X^{\prime}$ then we say that $f$ is $G$-homotopic to $g$ relative to $X^{\prime}$ if there exists a $G$ homotopy $\mathcal{H}\colon f\simeq_G g$ such that $\mathcal{H}\circ (i\times id) =\alpha \circ pr_1,$
where $\alpha=f|_{X^{\prime}}=g|_{X^{\prime}}$ and $pr_1 \colon X^{\prime}\times \Delta[1] \rightarrow X^{\prime}$ is the projection onto the first factor. In this case we write $\mathcal{H}\colon f\simeq_G g~(rel~ X^{\prime}).$
\end{Def}

\begin{Def}
A $G$-simplicial set is a $G$-Kan complex if for every subgroup $H\subseteq G$ the fixed point simplicial set $X^H$ is a Kan complex.
\end{Def}

\begin{Rem}\lb{equi}
Recall (\cite{ag}, \cite{fg}) that the category $G\mathcal{S}$ of $G$-simplicial sets and $G$-simplicial maps between $G$-simplicial sets has a closed model structure \cite{qui}, where the fibrant objects are the $G$-Kan complexes and cofibrant objects are the $G$-simplicial sets. From this it follows that $G$-homotopy on the set of $G$-simplicial maps $X\rightarrow Y$ is an equivalence relation, for every $G$-simplicial set $X$ and $G$-Kan complex $Y$. More generally, relative $G$-homotopy is an equivalence relation if the target is a $G$-Kan complex.
\end{Rem}

We consider $G/H \times \Delta[n]$ as a simplicial set where $(G/H \times \Delta[n])_q=G/H \times \Delta[n]_q$ with face and degeneracy maps as $id\times \partial_i$ and $id \times s_i$. Note that the group $G$ acts on $G/H$  by left translation. With this $G$-action  on the first factor and trivial action on the second factor $G/H \times \Delta[n]$ is a $G$-simplicial set.

Let $X$ be any $G$-simplicial set. A $G$-simplicial map $\sigma\colon G/H \times \Delta[n]\rightarrow X$ is called an equivariant $n$-simplex of type $H$ in $X$.
\begin{Rem}\lb{correspondence}
We remark that for a $G$-simplicial set $X$, the set of equivariant $n$-simplices of type $H$ in $X$ is in bijective correspondence with $n$-simplices of $X^H$. For an equivariant $n$-simplex $\sigma$, the corresponding $n$-simplex is $\sigma^{\prime}=\sigma (eH, \Delta_n).$ The simplicial map $\Delta[n]\rightarrow X^H,~~\Delta_n \mapsto \sigma^{\prime}$
will be denoted by $\overline{\sigma}.$
\end{Rem}

 We shall call $\sigma$  degenerate or non-degenerate according as the $n$-simplex $\sigma^{\prime} \in X^H_n$ is degenerate or non-degenerate.

Recall that the category of canonical orbits, denoted by $O_G,$ is a category whose  objects are cosets $G/H$, as $H$ runs over the all subgroups of $G$. A morphism from $G/H$ to $G/ K$ is a $G$-map. Such a morphism determines and is determined by a subconjugacy relation $a^{-1}Ha\subseteq K$ and is given by $\hat{a}(eH)=aK$. We denote this morphism by $\hat{a}$ \cite{br}.
\begin{Def}
A contravariant functor from $O_G$ to the category of simplicial sets $\mathcal{S}$ is called an $O_G$-simplicial set. A map between $O_G$-simplicial sets is a natural transformation of functors.
\end{Def}
We shall denote the category of $O_G$-simplicial sets by $O_G\mathcal{S}.$

For a commutative ring $\Lambda$, let $\Lambda$-alg denote the category of commutative $\Lambda$-algebras with unity and algebra homomorphisms preserving unity. The category of $\Lambda$-modules and module maps is denoted by $\Lambda$-mod. The category of chain complexes of $\Lambda$-modules is denoted by $ch_{\Lambda}.$ The notion of $O_G$-groups, $O_G$-$\Lambda$-algebras or $O_G$-chain complexes has the obvious meaning replacing $\mathcal{S}$ by $\mathcal{G}rp$ (the category of groups), $\Lambda$-alg or $ch_{\Lambda}$, respectively.

For any two $O_G$-simplicial sets (respectively, $O_G$-groups) $T$ and $T^{\prime}$, we define their product $(T\times T^{\prime})\in O_G\mathcal{S}$ (respectively, $O_G$-$\mathcal{G}rp$) as $$(T\times T^{\prime})(G/H)=T(G/H)\times T^{\prime}(G/H)$$ for objects $G/H$ of $O_G$ and $(T\times T^{\prime})(\hat{a})=T(\hat{a})\times T^{\prime}(\hat{a})$ for a morphism $\hat{a}$ of $O_G$.

For a $G$-simplicial set $X$, with a $G$-fixed $0$-simplex $v$, we have an $O_G$-group $\underline{\pi}X$ defined as follows. For any subgroup $H$ of $G$,
$$\underline{\pi}X(G/H):= \pi_1(X^H,v)$$ and for a morphism $\hat{a}\colon G/H\rightarrow G/K,~~a^{-1}Ha\subseteq K$, $\underline{\pi}X(\hat{a})$ is the homomorphism of fundamental groups induced by the simplicial map $a\colon X^K\rightarrow X^H.$

\begin{Def}
An $O_G$-group $\rho$ is said to act on an $O_G$-simplicial set ($O_G$-$\Lambda$-algebra or $O_G$-chain complex) $T$ if for every subgroup $H\subseteq G$, $\rho(G/H)$ acts on $T(G/H)$ and this action is natural with respect to maps of $O_G.$ Thus if
$$\phi (G/H)\colon \rho(G/H)\times T(G/H)\rightarrow T(G/H)$$ denotes the action of $\rho(G/H)$ on $T(G/H)$ then for each subconjugacy relation\\ $a^{-1}Ha\subseteq K,$
$$\phi (G/H)\circ (\rho(\hat{a})\times T(\hat{a}))= T(\hat{a})\circ \phi(G/K).$$
\end{Def}
\begin{Def}
 Let an $O_G$-group $\rho$ acts on $O_G$-simplicial sets $T$ and $T^{\prime}$. A map $f\colon T\rightarrow T^{\prime}$ is called $\rho$-equivariant if $$f(G/H)(ax)=af(G/H)(x),~a\in \rho(G/H),~x\in T(G/H),$$ for each subgroup $H$ of $G$.
\end{Def}

\begin{Def}
 Let $L,L^{\prime}$ be $O_G$-chain complexes. Two natural transformations $\mathsf{v}=\{\mathsf{v}_n\},\mathsf{w}=\{\mathsf{w}_n\}\colon L\rightarrow L^{\prime}$ are said to be homotopic if there exist natural transformations $$\mathcal{H}_n:\mathsf{v}_n\rightarrow \mathsf{w}_{n+1},~ n\geq 0$$ such that $\{\mathcal{H}_n(G/H)\}_{n\geq 0}$ is a chain homotopy of the chain maps $\mathsf{v}(G/H), \mathsf{w}(G/H)$ for each subgroup $H$ of $G$. Symbolically we write $\mathcal{H}:\mathsf{v}\simeq \mathsf{w}$.

If an $O_G$-group $\rho$ acts on $L,L^{\prime}$ and $\mathsf{v},\mathsf{w}$ are $\rho$-equivariant, then $\mathsf{v},\mathsf{w}$ are said to be $\rho$-equivariantly homotopic if there exists a homotopy $\mathcal{H}\colon \mathsf{v}\simeq \mathsf{w}$ which satisfies $$ \mathcal{H}_n(G/H)(ax)=a \mathcal{H}_n(G/H)(x)~\mbox{for}~a\in \rho(G/H),~x\in \mathsf{v}_n(G/H),~ H\subseteq G.$$
 \end{Def}
\begin{Def}
The tensor product $L\otimes L^{\prime}\colon O_G\rightarrow ch_{\Lambda}$ of two $O_G$-chain complexes $L$ and $L^{\prime}$ is defined as $$(L\otimes L^{\prime})(G/H)=L(G/H)\otimes L^{\prime}(G/H),$$ for each object $G/H$ of $O_G$ and $(L\otimes L^{\prime})(\hat{a})=L(\hat{a})\otimes L^{\prime}(\hat{a})$ for a morphism $\hat{a}$ of $O_G$.
\end{Def}
Note that a chain complex $W$ can be considered as an $O_G$-chain complex in the trivial way, that is, $W(G/H)=W,~W(\hat{a})=id.$ So the tensor product of $W$ with an $O_G$-chain complex is defined.

Throughout the paper, unless otherwise mentioned explicitly, all the tensor products are over the ring $\Lambda.$
\section{Simplicial Bredon-Illman cohomology with local\\ coefficients}

In this section we recall \cite{ms} the relevant notion of fundamental groupoid of a $G$-simplicial set $X$, the notion of equivariant local coefficients on $X$ and
the definition of simplicial Bredon-Illman cohomology with local coefficients.

We begin with the notion of fundamental groupoid.  Recall \cite{gj} that the fundamental groupoid $\pi X$ of a Kan complex $X$ is a category having as objects all $0$-simplexes of $X$ and a morphism $ x\rightarrow y$ in $\pi X$ is a homotopy class of $1$-simplices $\omega \colon  \Delta [1] \rightarrow X$ rel $\partial \Delta [1]$ such that $\omega \circ \delta_0 = \overline{y}$, $\omega \circ \delta_1 = \overline{x}$. If $\omega_2$ represents an arrow from $x$ to $y$ and $\omega_0$ represents an arrow from $y$ to $z$,  then their composite $[\omega_0]\circ [\omega_2]$ is represented by $\Omega \circ \delta_1$, where the simplicial map $\Omega\colon \Delta[2]\rightarrow X$ corresponds to a $2$-simplex, which is determined by the compatible pair $(\omega^{\prime}_0,~~, \omega^{\prime}_2)$. For a simplicial set $X$ the notion of fundamental groupoid is defined via the geometric realization and the total singular functor.

The fundamental groupoid of a $G$-simplicial set is defined as follows.
\begin{Def}
Let $X$ be a $G$-Kan complex. The fundamental groupoid $\Pi X$ is a category with objects equivariant $0$-simplices
$$x_{H}\colon G/ H \times \Delta[0] \rightarrow X$$
of type $H$, as $H$ varies over all subgroups of $G$. Given two objects $x_H$ and  $y_K$ in $\Pi X$, a morphism from $x_H \longrightarrow y_K$ is defined as follows. Consider the set of all pairs $(\hat{a},\phi)$ where $\hat{a}\colon G/H \rightarrow G/K$ is a morphism in $O_G$, given by a subconjugacy relation $a^{-1}Ha\subseteq K$, $a\in G$ so that $\hat{a}(eH)=aK$ and $\phi\colon G/H \times \Delta[1]\rightarrow X$ is an equivariant $1$-simplex such that
$$\phi \circ (id\times \delta_1) = x_H,~~\phi \circ (id\times \delta_0)=y_K\circ (\hat{a}\times id).$$

The set of morphisms in $\Pi X$ from $x_H$ to $y_K$ is a quotient of the set of pairs mentioned above by an equivalence relation $` \sim\mbox{'},$ where $(\hat{a}_{1},\phi_{1})\sim(\hat{a}_{2},\phi_{2})$ if and only if  $a_1=a_2=a$ (say) and there exists a $G$-homotopy
$\mathcal{H} \colon  G/H \times \Delta [1] \times \Delta [1] \rightarrow X$ of $G$-homotopies such that $\mathcal{H} \colon  \phi_1 \simeq_G \phi_2$ (rel $G/H \times \partial \Delta [1]$). Since $X$ is a $G$-Kan complex, by Remark \ref{equi}, $\sim$ is an equivalence relation. We denote the equivalence class of $(\hat{a},\phi)$ by $[\hat{a},\phi]$. The set of equivalence classes is the set of morphisms in $\Pi X$ from $x_H$ to $y_K$.

The composition of morphisms in $\Pi X$ is defined as follows. Given two morphisms
$$
\begin{diagram}
\node{x_{H}} \arrow{e,t}{[\hat{a}_{1},\phi_{1}]} \node{y_{K}} \arrow{e,t}{[\hat{a}_{2},\phi_{2}]} \node{z_{L}}
\end{diagram}
$$
their composition $[\hat{a}_2, \phi_2]\circ [\hat{a}_1, \phi_1]$ is $[\widehat{a_1a_2}, \psi]\colon  x_H \rightarrow z_L$, where the first factor is the composition
\[
\begin{diagram}
\node{G/ H} \arrow[2]{e,t}{\hat{a}_{1}} \node[2]{G/ K} \arrow[2]{e,t}{\hat{a}_{2}} \node[2]{G/L}
\end{diagram}
\]
and $\psi\colon  G/H \times \Delta[1] \rightarrow X$ is an equivariant $1$-simplex of type $H$  as described below. Let $x$ be a $2$-simplex in the Kan complex $X^H$ determined by the compatible pair of $1$-simplices $(a_1\phi_2^{\prime}, ~~, \phi_1^{\prime})$ so that $ \partial_0x = a_1\phi_2^{\prime}$ and $\partial_2x= \phi_1^{\prime}$. Then $\psi$ is given by $\psi (eH, \Delta_1) = \partial_1x$.
\end{Def}
Observe that $\phi^{\prime}$ is a $1$-simplex in $X^H$  such that $\partial_1\phi^{\prime} = x_H^{\prime}$ and $\partial_0\phi^{\prime} = ay_K^{\prime}$. Moreover the $0$-simplex $ay_K^{\prime}$ in $X^H$ corresponds to the composition
$$G/H \times \Delta [0]\xrightarrow{\hat{a}\times id} G/K \times \Delta [0]\xrightarrow{y_K} X$$ and $\phi $ is a $G$-homotopy $x_H \simeq_G y_K \circ (\hat{a}\times id)$ (cf. Remark \ref{correspondence} for notations).

It is proved in \cite{ms} that the composition is well defined. For a version of fundamental groupoid of a $G$-space we refer \cite{luck}, \cite{mm}.

Observe that if $X$ is a $G$-simplicial set then $S|X|$ is a $G$-Kan complex, where for any space $Y$, $SY$  denotes the total singular complex and for any simplicial set $X$, $|X|$ denotes the geometric realization of $X$.
\begin{Def}
For a $G$-simplicial set $X$, we define the fundamental groupoid $\Pi X$ of $X$ by $\Pi X := \Pi S|X|.$
\end{Def}
Note that if $F\colon X\rightarrow Y$ is a $G$-simplicial map then there exists an obvious induced functor $\Pi(F)\colon \Pi X\rightarrow \Pi Y$ which assigns to each object $x_H$ of $\Pi X$, the object $F\circ x_H$ of $\Pi Y$ and a morphism $[\hat{a},\phi]$ in $\Pi X$ to the morphism $[\hat{a},F\circ \phi]$ of $\Pi Y$.

\begin{Rem}\lb{morphism}
 Suppose $\xi$ is a morphism from $x$ to $y$ in $\pi X^H$, given by a homotopy class $[\overline{\omega}],$ where $\overline{\omega}\colon \Delta[1]\lgra X^H$ represents the $1$-simplex in $X^H$ from $x$ to $y$. Let $x_H$ and $y_H$ be the objects in $\pi X^H$ defined respectively by $$x_H(eH,\Delta_0)=x,~~y_H(eH,\Delta_0)=y.$$ Then we have a morphism $[id, \omega]\colon x_H\lgra y_H$ in $\Pi X$, where $\omega(eH,\Delta_1) = \overline{\omega}(\Delta_1).$  We shall denote this morphism corresponding to $\xi$ by $b\xi.$
\end{Rem}

\begin{Def}
An equivariant local coefficients on a $G$-simplicial set $X$ is a contravariant functor from $\Pi X$ to the category $\Lambda$-alg.
\end{Def}

Next, we briefly describe  the simplicial version of Bredon-Illman cohomology with local coefficients as introduced in \cite{ms}.

Let $X$ be a $G$-simplicial set and $M$ an equivariant local coefficients on $X$. For each equivariant $n$-simplex $\sigma\colon G/H\times \Delta[n]\rightarrow X,$ we associate an equivariant $0$-simplex $\sigma _{H}\colon G/H\times\Delta[0]\rightarrow X$ given by
$$\sigma_H= \sigma\circ (id\times \delta_{(1,2,\cdots,n)}),$$
where $\delta_{(1,2, \cdots ,n)}$ is the composition
$$\delta_{(1,2, \cdots ,n)}\colon  \Delta[0]\xrightarrow{\delta_1} \Delta[1]\stackrel{\delta_2}{\rightarrow}\cdots \stackrel{\delta_ n}{\rightarrow} \Delta[n].$$
The $j$-th face of $\sigma$ is an equivariant $(n-1)$-simplex of type $H$, denoted by $\sigma^{(j)}$, and is defined by
$$ \sigma^{(j)}= \sigma \circ (id \times \delta_j), ~0\leq j\leq n.$$
\begin{Rem}\lb{initial}
Note that $\sigma^{(j)}_{H}=\sigma_{H}\mbox{ for }j> 0,$ and
$$\sigma^{(0)}_H = \sigma \circ (id \times \delta_{(0,2,\cdots, n)}).$$
\end{Rem}

Let $C^{n}_{G}(X;M)$ be the $\Lambda$-module of all functions $f$ defined on equivariant $n$-simplexes $\sigma: G/H\times\Delta[n]\rightarrow X$ such that $f(\sigma)\in M(\sigma _{H})$ with $f(\sigma)=0,$ if $\sigma$ is degenerate. We have a morphism $\sigma_*=[id,\alpha]$ in $\Pi X$ from $\sigma_{H}$ to $\sigma^{(0)}_{H}$ induced by $\sigma$, where $\a \colon  G/H \times \Delta[1] \rightarrow X$ is given by $\a = \sigma \circ (id \times \delta_{(2,\cdots ,n)}).$ Define a homomorphism $$\delta\colon C_{G}^{n}(X;M)\rightarrow C_{G}^{n+1}(X;M);~f\mapsto \delta f$$
where for any equivariant $(n+1)$-simplex $\sigma $ of type $H$,
$$(\delta f)(\sigma)=M(\sigma_{*})(f(\sigma^{(0)}))+\sum_{j=1}^{n+1}(-1)^{j}f(\sigma^{(j)}).$$
A routine verification shows that $\delta\circ \delta=0.$ Thus $\{C_{G}^{*}(X;M),\delta \}$ is a cochain complex. We are interested in a subcomplex of this cochain complex as described below.

Let $\eta\colon G/H\times \Delta[n]\rightarrow X$ and $\tau\colon G/ K\times \Delta[n]\rightarrow X$ be two equivariant $n$-simplexes. Suppose there exists a $G$-map $\hat{a} \colon G/H\rightarrow G/K,~ a^{-1}Ha\subseteq K,$ such that $\tau \circ (\hat{a}\times id ) = \eta.$ Then $\eta$ and $\tau $ are said to be compatible under $\hat{a}$. Observe that if $\eta$ and $\tau$ are compatible as described above then $\eta$ is degenerate if and only if $\tau$ is degenerate. Moreover notice that in this case, we have a morphism
$[\hat{a},k]\colon \eta_{H}\rightarrow \tau_{K}$ in $\Pi X$, where $k = \eta_H \circ (id \times \sigma_0),$ where $\sigma_0\colon \Delta[1]\lgra \Delta[0]$ is the simplicial map as described in Section 2. Let us denote this induced morphism by $a_*$.

\begin{Def}
We define $S_{G}^{n}(X;M)$ to be the submodule of $C_{G}^{n}(X;M)$ consisting of all functions f such that if $\eta$ and $\tau$ are equivariant n-simplexes in X which are compatible under $\hat{a} $ then $f(\eta)=M(a_*)(f(\tau))$.
\end{Def}

If $f\in S_{G}^{n}(X;M)$ then one can verify that $\delta f \in S_{G}^{n+1}(X;M).$ Thus we have a cochain complex of $\Lambda$-modules $S_G(X;M) = \{S^n_G(X;M), \delta\}.$

\begin{Def}
Let X be a G-simplicial set with equivariant local coefficients M on it. Then the $n$-th Bredon-Illman cohomology of $X$ with local coefficients $M$ is defined by $$H^n_G(X;M):=H^{n}(S_{G}(X;M)).$$
\end{Def}

Suppose that $X,Y$ are $G$-simplicial sets and $M,N$ are equivariant local coefficients on $X$ and $Y$ respectively. A map from $(X,M)$ to $(Y,N)$ is a pair $(F,\gamma)$, where $F\colon X\rightarrow Y$ is a $G$-simplicial map and $\gamma\colon N\circ\Pi(F)\rightarrow M$ is a natural transformation of functors, $\Pi(F)\colon \Pi X\rightarrow \Pi Y$ being the map induced by $F$. A map $(F,\gamma)\colon (X,M)\rightarrow (Y,N)$ naturally induces a cochain map $(F,\gamma)^{\#}\colon S_G^*(Y;N)\rightarrow S_G^*(X;M)$ as follows. For $f\in S_G^*(Y;N)$ and an equivariant $n$-simplex $\sigma$ in $X$ of type $H$, $(F,\gamma)^{\#}(f)(\sigma)=\gamma(\sigma_H)f(F\circ\sigma)$. Therefore we have an induced map $(F,\gamma)^*\colon H_G^*(Y;N)\rightarrow H_G^*(X;M)$ in cohomology.

We now define the cup product in simplicial Bredon-Illman cohomology with local coefficients. Let $\sigma\colon G/H\times \Delta[n+m]\rightarrow X$ be an equivariant $(n+m)$-simplex of type $H$. Then define $\sigma\rfloor_n=\sigma\circ(id_{G/H}\times \delta_{(n+1,\cdots,n+m)})$, $\lfloor_m\sigma=\sigma\circ (id_{G/H}\times \delta_{(0,\cdots,n)})$ where $\delta_{(n+1,\cdots,n+m)}\colon \Delta[n]\rightarrow \Delta[n+m]$ and $\delta_{(0,\cdots,n)}\colon \Delta[m]\rightarrow \Delta[n+m]$ are defined as before. For cochains $f\in S_G^n(X;M)$ and $g\in S_G^m(X;M)$, the cup product $f\cup g \in S_G^{n+m}(X;M)$ is the cochain whose value on $\sigma$ is given by the formula $$(f\cup g)(\sigma)=f(\sigma\rfloor_n)(M(\sigma_{n+1})g(\lfloor_m\sigma))$$ where $\sigma_{n+1}=[id,\sigma\circ(id_{G/H}\times \delta_{(1,\cdots,n,n+2,\cdots, n+m)})]$ is a  morphism in $\Pi X$ from $(\sigma\rfloor_n)_H$ to $(\lfloor_m\sigma)_H.$ A routine verification shows that $f\cup g$ belongs to $S_G^{n+m}(X;M)$ and $$d(f\cup g)=df\cup g+(-1)^{deg(f)}f\cup dg.$$ Therefore it induces a product in cohomology which is associative and graded commutative. Thus $H^*_G(X;M)$ is an associative graded algebra.

Suppose $M$ is an equivariant local coefficients on a $G$-simplicial set $X$ with a $G$-fixed $0$-simplex $v$. Then $M$ determines an $O_G$-$\Lambda$-algebra $M_0$ equipped with an action of the $O_G$-group $\underline{\pi}X$ as described below.

For any subgroup $H$ of $G$, let $v_{H}$ be the object of type $H$ in $\Pi X$ defined by $$v_{H}\colon G/H\times \Delta[0]\rightarrow X,~v_H(eH,\Delta_0)=v.$$
Then for any morphism $\hat{a}\colon G/H\rightarrow G/ K$ in $O_{G}$ given by a subconjugacy relation $a^{-1}Ha\subseteq K$, we have a morphism  $[\hat{a},k]\colon v_{H}\rightarrow v_{K}$ in $\Pi X,$ where $k \colon G/H\times \Delta[1] \rightarrow X$ is given by $k(eH, \Delta_1)= s_0v$.
Define an $O_G$-$\Lambda$-algebra $M_{0}$ by $$M_0(G/H):=M(v_{H}),~H\subseteq G,$$ and $M_{0}(\hat{a})=M[\widehat{a},k]$ for a morphism $\hat{a}$ in $O_G$.

We now describe the action of the $O_G$-group $\underline{\pi}X$ on $M_0$. Let $\alpha=[\overline{\phi}]\in \underline{\pi}X(G/H)=\pi_{1}(X^{H},v)$. Then the morphism $[id,\phi]\colon v_{H}\rightarrow v_{H}$, determined by $\phi(eH,\Delta_1)= \overline{\phi}(\Delta_1),$ is an equivalence in the category $\Pi X$. This yields a group homomorphism
$$b\colon \pi_{1}(X^{H},v)\rightarrow Aut_{\Pi X}(v_{H}), ~\a=[\overline{\phi}] \mapsto b(\a)=[id,\phi].$$ The composition of the map $b$ with the group homomorphism $Aut_{\Pi X}(v_{H})\rightarrow Aut_{\Lambda\mbox{-alg}}(M(v_{H}))$, which sends $\a\in Aut_{\Pi X}(v_{H})$ to $[M(\a)]^{-1}$, defines the action of $\pi_{1}(X^{H},v)$ on $M_{0}(G/H)$. It is routine to check that this action is natural with respect to morphisms of $O_G$.

Conversely, an $O_G$-$\Lambda$-algebra $M_0$, equipped with an action of the $O_G$-group $\underline{\pi}X,$ defines an equivariant local coefficients $M$ on $X$, where $X$ is $G$-connected and $v\in X^G$ a fixed $0$-simplex \cite{ms}.

\section{Eilenberg Theorem}
In this section we recall a version of Eilenberg theorem \cite{ms1} for simplicial Bredon-Illman cohomology with local coefficients.

Let $\mathcal{A}_{\Lambda}$ denote the category with objects the triples $(T,M_0,\rho)$ where $T$ is an $O_G$-simplicial set, $M_0$ an $O_G$-$\Lambda$-algebra and $\rho$ is an $O_G$-group which operates on both $T$ and $M_0$. A morphism from $(T,M_0,\rho)$ to $(T^{\prime},M_0^{\prime},\rho^{\prime})$ is a triple $(f_0,f_1,f_2)$, where $f_0\colon T\rightarrow T^{\prime}$, $f_1\colon M_0^{\prime}\rightarrow M_0$ and $f_2\colon \rho\rightarrow \rho^{\prime}$ are maps in the appropriate categories such that $$f_0(G/H)(\alpha x)=f_2(G/H)(\alpha)f_0(G/H)(x),~f_1(G/H)[f_2(G/H)(\alpha)m_0^{\prime}]=\alpha f_1(G/H)(m_0^{\prime}),$$  $H\subseteq G,~ x\in T(G/H),\alpha\in \rho(G/H),m_0^{\prime}\in M_0^{\prime}(G/H).$

 The $\rho$-equivariant cohomology of $T$ with coefficients $M_0$ is defined as follows. We have an $O_G$-chain complex $\{\underline{C}_*(T),\partial_*\}$, defined by $$\underline{C}_n(T)\colon O_G\rightarrow \Lambda\mbox{-mod},~~G/H\mapsto C_n(T(G/H);\Lambda),$$
where $C_n(T(G/H);\Lambda)$ is the free $\Lambda$-module generated by the non-degenerate $n$-simplices of $T(G/H)$. For any morphism $\hat{a}\colon G/H\rightarrow G/K$ in $O_G$, $$\underline{C}_n(T)(\hat{a})=a_{\#}\colon C_n(T(G/K);\Lambda)\rightarrow C_n(T(G/H);\Lambda)$$ is induced by the simplicial map $T(\hat{a})\colon T(G/K)\rightarrow T(G/H).$ The boundary $\partial_n\colon \underline{C}_n(T)\rightarrow \underline{C}_{n-1}(T)$ is a natural transformation defined by $\partial_n(G/H)\colon C_n(T(G/H);\Lambda)\rightarrow C_{n-1}(T(G/H);\Lambda),$ where $\partial_n(G/H)$ is the ordinary boundary map of the simplicial set $T(G/H)$. The action of $\rho$ on $T$ induces an action of $\rho$ on the $O_G$-chain complex $\{\underline{C}_*(T),\partial_*\}$. We form the cochain complex $$\{C^*_{\rho}(T;M_0)=Hom_{\rho}(\underline{C}_*(T),M_0),\delta^* \}$$ where $Hom_{\rho}(\underline{C}_n(T),M_0)$ consists of all natural transformations $\underline{C}_n(T)\xrightarrow{f} M_0$ respecting the action of $\rho$ and $\delta^n f$ is given by $f\circ \partial_{n+1}.$ Then the $n$-th $\rho$-equivariant cohomology of $T$ with coefficients $M_0$ is given by $$H^n_{\rho}(X;M_0):= H_n(C^*_{\rho}(T;M_0)).$$
\begin{Rem}\lb{map}
 It is easy to observe that a morphism $(f_0,f_1,f_2)\colon (T,M_0,\rho)\rightarrow (T^{\prime},M_0^{\prime},\rho^{\prime})$ induces a cochain map $C^*(f_0,f_1,f_2)\colon C^*_{\rho}(T;M_0)\rightarrow C^*_{\rho^{\prime}}(T^{\prime};M_0^{\prime}).$
\end{Rem}

The cochain complex $C^*_{\rho}(T;M_0)$ is equipped with a cup product, defined as follows. We have a natural transformation $$\underline{\xi}\colon \underline{C}_*(T\times T)\rightarrow \underline{C}_*(T)\otimes\underline{C}_*(T),$$ where $\underline{\xi}(G/H)$ is the Alexander-Whitney map for the simplicial set $T(G/H)$, $H\subseteq G$ \cite{may}. We have a $\rho$-action on $\underline{C}(T)$, induced by the $\rho$-action on $T$ and hence diagonal actions of $\rho$ on $T\times T$ and $\underline{C}(T)\otimes \underline{C}(T)$. Since the Alexander-Whitney map of simplicial sets is a natural map, $\underline{\xi}$ is equivariant with the induced actions of $\rho$ on $\underline{C}(T\times T)$ and $\underline{C}(T)\otimes \underline{C}(T)$. Then the cup product is defined as the composition
 $$C_{\rho}^*(T;M_0)\otimes C_{\rho}^*(T;M_0)\xrightarrow{\alpha} Hom_{\rho}(\underline{C}(T)\otimes \underline{C}(T),M_0)\xrightarrow{\underline{\xi}^*} C_{\rho}^*(T\times T;M_0)\xrightarrow{D^*} C_{\rho}^*(T;M_0).$$
Here $\alpha\colon  C_{\rho}^*(T;M_0)\otimes C_{\rho}^*(T;M_0)\rightarrow Hom_{\rho}(\underline{C}_*(T\times T),M_0)$ is defined by $$\alpha(f\otimes g)(G/H)(x\otimes y)=(-1)^{deg(x)deg(y)}f(G/H)(x)g(G/H)(y)$$ where $f,g\in C_{\rho}^*(T)$ and $x,y\in \underline{C}_*(T)(G/H)$ and $D\colon T\rightarrow T\times T$ is the diagonal map.
\begin{Rem}\lb{4.1}
 The cochain complex $C_{\rho}^*(T;M_0)$, equipped with the above cup product, is a homotopy associative differential $\Lambda$-algebra and the induced product in the cohomology is associative and graded commutative.
\end{Rem}

We now relate the simplicial Bredon-Illman cohomology with local coefficients of a one vertex $G$-Kan complex with the equivariant cohomology of its universal $O_G$-covering complex \cite{ms1}.

Let $X$ be a one vertex $G$-Kan complex. We denote the $G$-fixed vertex by $v$. Let $M$ be an equivariant local coefficients on $X$ and $M_0$ be the associated $O_G$-$\Lambda$-algebra as described at the end of the previous section. For any subgroup $H$ of $G$, let $$p_H\colon \widetilde{X^H}\rightarrow X^H$$ be the universal cover \cite{gugg,ms1} of $X^H$. The left translation $a\colon X^K\rightarrow X^H,$ corresponding to a G-map $\hat{a}\colon G/H\rightarrow G/ K$, $a^{-1}Ha\subseteq K,$ induces a simplicial map $\tilde{a}\colon \widetilde{X^K}\rightarrow \widetilde{X^H}$ such that $p_H\circ   \tilde{a}=a\circ p_K$. This defines an $O_G$-Kan complex $\widetilde{X}$ by setting $\widetilde{X}(G/H)=\widetilde{X^H}$ and $\widetilde{X}(\hat{a})=\tilde{a}$. This is called the universal $O_G$-covering complex of $X$. This is the simplicial analogue of $O_G$ covering space as introduced in \cite{mm}. We refer \cite{luck} for a more general version, called `universal covering functor'. The natural actions of $\underline{\pi}X(G/H) =\pi_1(X^H, v)$ on $\widetilde{X}(G/H)= \widetilde{X^H}$ as $H$ varies over subgroups of $G$, define an action of the $O_G$-group $\underline{\pi}X$ on $\widetilde{X}$.
Thus $(\tilde{X},M_0,\underline{\pi}X)$ is an object of $\mathcal{A}_{\Lambda}.$

\begin{Thm}\cite{ms1}\lb{eilenberg}
 Let $X$ be a one vertex $G$-Kan complex with equivariant local coefficients $M$ on it. Then, with notations as above, there exists an isomorphism of graded algebras
$$H_G^*(X;M)\cong H_{\underline{\pi}X}^*(\tilde{X};M_0).$$
\end{Thm}

The proof is obtained by constructing isomorphism at the cochain level.
The explicit isomorphism is described as follows \cite{ms1}. Define
$$\mu \colon  S^n_G(X;M) \rightarrow Hom_{\underline{\pi}X}(\underline{C}_n(\widetilde{X}), M_0)$$ as follows. Let $f\in S^n_G(X; M)$ and $y$ be a non-degenerate $n$-simplex in $\widetilde{X^H}.$ Let $\sigma $ be the equivariant $n$-simplex of type $H$ in $X$ such that $\overline{\sigma} = p_H\circ \overline{y},$ where $\overline{y} \colon \Delta[n] \rightarrow \widetilde{X^H}$ is the simplicial map with $\overline{y}(\Delta_n)=y.$ Then $\mu(f) \in Hom_{\underline{\pi}X}(\underline{C}_n(\widetilde{X}), M_0)$ is given by
$$\mu (f)(G/H)(y) = M(b\xi_H(\partial_{(1,2,\cdots ,n)}y))f(\sigma),$$ where $\partial_{(1,2,\cdots ,n)}y=\partial_1\partial_2\cdots\partial_n y.$

The inverse of $\mu$, $$\mu^{-1}\colon Hom_{\underline{\pi}X}(\underline{C}_n(\tilde{X}),M_0) \rightarrow  C_G^n(X;M)$$ is described as follows. Let $f\in Hom_{\underline{\pi}X}(\underline{C}_n(\tilde{X}),M_0)$ and $\sigma$ be a non-degenerate equivariant $n$-simplex of type $H$ in X. Choose an $n$-simplex $y$ in $\widetilde{X^H}$ such that $p_H(y) = \sigma (eH, \Delta_n).$ Then $\mu^{-1}(f)$ is given by
$$\mu^{-1} (f)(\sigma ) = M(b\xi_H(\partial_{(1,2, \cdots, n)}y))^{-1}f(G/H)(y).$$
It is easy to check that $\mu(f\cup g)=\mu(f)\cup \mu(g)$ for $f,g\in S_G^*(X;M).$ Hence we have an isomorphism $$\mu^*:H_G^*(X;M)\cong H_{\underline{\pi}X}^*(\tilde{X};M_0)$$ of graded $\Lambda$-algebras.

\section {Steenrod reduced power operations}
In this section we briefly recall the relevant part of the general algebraic approach to Steenrod operations by P. May \cite{may1}, necessary for our purpose. We apply this method to construct Steenrod power operations in equivariant cohomology of $O_G$-simplicial sets in general. In particular, for a one vertex $G$-Kan complex $X$, we have reduced power operations defined for $\underline{\pi}X$-equivariant cohomology of the universal $O_G$-covering complex $\tilde{X}$. We then apply the Theorem \ref{eilenberg} to deduce the Steenrod power operations in simplicial Bredon-Illman cohomology with local coefficients.

 Let $p$ be an odd prime and $\Lambda$ be the commutative ring $\mathbb{Z}_p.$ Suppose $\Sigma_p$ denotes the symmetric group on $p$-letters and $\pi$ be a subgroup of $\Sigma_p$. Let $\Lambda\pi$ (respectively, $\Lambda\Sigma_p$) denote the group ring of $\pi$ (respectively, $\Sigma_p$) over $\Lambda$ and $V,W$ be free resolutions of $\Lambda$ over $\Lambda\Sigma_p,\Lambda\pi$ respectively. If $\pi$ is cyclic of order $p$ with generator $\alpha=(p,1,2,\cdots,p-1)$, we shall use the following canonical model of $W.$
 Let $W_i$ be $\Lambda\pi$-free module on one generator $e_i,i\geq 0.$  Let $N=1+\alpha+\cdots +\alpha^{p-1}$ and $T=\alpha -1$ in $\Lambda\pi.$ Define differential $d$, augmentation $\epsilon\colon W_0\rightarrow \Lambda$, and coproduct $\psi$ on $W$ by the formulas
$$d(e_{2i+1})=Te_{2i},~d(e_{2i})=Ne_{2i-1},~ \epsilon(\alpha^je_0)=1;$$
$$\psi(e_{2i+1})=\sum_{j+k=i}e_{2j}\otimes e_{2k+1}+ \sum_{j+k=i}e_{2j+1}\otimes \alpha e_{2k},$$
$$\psi(e_{2i})=\sum_{j+k=i}e_{2j}\otimes e_{2k}+ \sum_{j+k=(i-1)}\sum_{0\leq r<s<p}\alpha^r e_{2j+1}\otimes \alpha^s e_{2k}.$$
Thus $W$ is a differential $\Lambda\pi$-coalgebra and a $\Lambda\pi$-free resolution of $\Lambda.$

From now on, unless otherwise stated, $\pi$ will be cyclic of order $p$ with generator $\alpha.$

By a $\Lambda$-complex K, we will mean a $\mathbb{Z}$-graded cochain complex of $\Lambda$-modules with differential of degree $1$. We denote the $p$-fold tensor product $K\otimes \cdots\otimes K$ by $K^p$. Then $K^p$ becomes a $\Lambda\pi$-complex by the following $\pi$ operation,
$$\tau(u_1\otimes\cdots \otimes u_p)=\gamma(\tau)u_1\otimes\cdots u_{i-1}\otimes u_{i+1}\otimes u_{i}\otimes u_{i+2}\cdots \otimes u_p$$ where $\gamma(\tau)=(-1)^{deg(u_i)deg(u_{i+1})}$ if $\tau$ is the interchange of $i$-th and $(i+1)$-th factor. We consider $W$ as a non-positively graded $\Lambda$-complex. The inclusion of $\pi$ in $\Sigma_p$ induces a morphism $j\colon W\rightarrow V$ of $\Lambda\pi$-complexes.

We have the following algebraic category $\mathfrak{C}(p)$ on which the Steenrod operations are defined. The objects of this category are pairs $(K,\theta),$ where $K=\{K^i\}_{i\in \mathbb{Z}}$ is a homotopy associative differential $\Lambda$-algebra, and $\theta\colon W\otimes K^p\rightarrow K$ is a morphism of $\Lambda\pi$-complexes, satisfying the following two conditions.

\begin{enumerate}
 \item The restriction of $\theta$ to $e_0\otimes K^p$ is $\Lambda$-homotopic to the iterated product $K^p\rightarrow K$, associative in some order.
 \item The morphism $\theta$ is $\Lambda\pi$-homotopic to a composite $W\otimes K^p\xrightarrow{j\otimes 1} V\otimes K^p\xrightarrow{\phi}K$ where $\phi$ is a morphism of $\Lambda\Sigma_p$-complexes.
\end{enumerate}

A morphism $f\colon (K,\theta)\rightarrow (K^{\prime},\theta^{\prime})$ is a morphism of $\Lambda$-complexes $f\colon K\rightarrow K^{\prime}$ such that the following diagram is $\Lambda\pi$-homotopy commutative.
$$
\begin{diagram}
\node{W\otimes K^p} \arrow{s,l}{id\otimes
f^p}\arrow{e,t}{\theta}\node{K}\arrow{s,r}{f} \\
\node{W\otimes (K^{\prime})^p}\arrow{e,t}{\theta^{\prime}}\node{K^{\prime}}
\end{diagram}
$$

The tensor product of two objects $(K,\theta)$ and $(K^{\prime},\theta^{\prime})$ is the pair $(K\otimes K^{\prime},\tilde{\theta})$, where $\tilde{\theta}$ is the composite
$$W\otimes (K\otimes K^{\prime})^p\xrightarrow{\psi\otimes \tilde{U}} W\otimes W \otimes K^p\otimes K^{\prime p}\xrightarrow{id\otimes \tilde{t}\otimes id} W\otimes K^p\otimes W\otimes K^{\prime p}\xrightarrow{\theta\otimes \theta^{\prime}} K\otimes K^{\prime}.$$
Here $\psi\colon W\rightarrow W\otimes W$ is the coproduct, $\tilde{U}\colon (K\otimes K^{\prime})^p\rightarrow K^p\otimes K^{\prime p}$ is the shuffling isomorphism and $\tilde{t}(x\otimes y)=(-1)^{deg(x)deg(y)} y\otimes x.$
\begin{Def}
 An object $(K,\theta)\in\mathfrak{C}(p)$ is said to be a Cartan object if the product $K\otimes K\rightarrow K$ is a morphism from $(K\otimes K,\tilde{\theta})$ to $(K,\theta).$
\end{Def}

For an object $(K,\theta)$ of $\mathfrak{C}(p)$, there are maps $D_i\colon H^q(K)\rightarrow H^{pq-i}(K), i\geq 0,$ defined as follows. For $x\in H^q(K)$, $e_i\otimes x^p$ is an well defined element of $H^{pq-i}(W\otimes_{\Lambda\pi}K^p)$ \cite{may1} and define $D_i(x)=\theta_*(e_i\otimes x^p),$ where $\theta_*\colon H^{pq-i}(W\otimes_{\Lambda\pi}K^p)\rightarrow H^{pq-i}(K)$ is induced by $\theta.$ We make the convention that $D_i=0$ for $i<0.$ Then the Steenrod power operations $$\mathcal{P}^s\colon H^q(K)\rightarrow H^{q+2s(p-1)}(K),~~\beta\mathcal{P}^s\colon H^q(K)\rightarrow H^{q+2s(p-1)+1}(K)$$ are defined by the following formulas
$$\mathcal{P}^s(x)=(-1)^r(m!)^q D_{(q-2s)(p-1)}(x),~~\beta\mathcal{P}^s(x)=(-1)^r(m!)^q D_{(q-2s)(p-1)-1}(x)$$ where $m=(p-1)/2$ and $r=s+m(q+q^2)/2.$

\begin{Prop}
 The power operations satisfies the following properties.
\begin{enumerate}
 \item $\mathcal{P}^s$ and $\beta\mathcal{P}^s$ are natural homomorphisms.
\item  $\mathcal{P}^s(x)=0$ if $2s>q$, $\beta\mathcal{P}^s=0$ if $2s\geq q$, and $\mathcal{P}^s(x)=x^p$ if $2s=q.$
\item If $(K,\theta)$ is a Cartan object then $\mathcal{P}^s$ satisfies the Cartan formulas
$$\mathcal{P}^s(xy)=\sum_{i+j=s}\mathcal{P}^i(x)\mathcal{P}^j(y),$$
$$\beta\mathcal{P}^{s+1}(xy)=\sum_{i+j=s}[\beta\mathcal{P}^{i+1}(x)\mathcal{P}^j(y)+(-1)^{deg(x)}\mathcal{P}^i(x)\beta\mathcal{P}^{j+1}(y)].$$
\end{enumerate}
\end{Prop}

\begin{Rem}
In general $\beta\mathcal{P}^s$ is single notation. But if $(K,\theta)$ is reduced mod $p$ (\cite{may1}) then Bockstein homomorphism $$\beta\colon H^n(K)\rightarrow H^{n+1}(K),$$can be defined and $\beta\mathcal{P}^s$ is the composition of $\mathcal{P}^s$ with the Bockstein.
\end{Rem}

Next we recall the definition of `Adem object' in $\mathfrak{C}(p)$ \cite{may1}. We need the following notations for the definition.

Consider $\Sigma_{p^2}$ as permutations on the $p^2$ symbols $\{(i,j)|1\leq i,j\leq p\}.$  Embed $\pi=<\alpha>$ ($\subseteq\Sigma_p$) in $\Sigma_{p^2}$ by letting $\alpha(i,j)=(i,j+1)$. Let $\alpha_i\in \Sigma_{p^2},1\leq i\leq p$ be defined as $\alpha_i(i,j)=(i,j+1)$ and $\alpha_i(k,j)=(k,j)$ for $k\neq i.$ Let $$\beta= \alpha_1\cdots \alpha_p,~ \nu=<\beta>,~\sigma=\pi\nu,~ \tau=<\alpha_1,\cdots,\alpha_p,\alpha>.$$ Note that $\beta$ and $\alpha_i$ are of order $p$ and the following relations hold. $$\alpha\alpha_i=\alpha_{i+1}\alpha;~\alpha_i\alpha_j=\alpha_j\alpha_i;~\alpha\beta=\beta\alpha.$$

Let $W_1=W$ and $W_2=W$ regarded as $\Lambda\pi$-free and $\Lambda\nu$-free resolutions of $\Lambda$ respectively. Let $\nu,\pi$
 operate trivially on $W_1,W_2$ respectively. Then $W_1\otimes W_2$ is a $\Lambda\sigma$-free resolution of $\Lambda$ with the diagonal action of $\sigma$ on $W_1\otimes W_2.$

If $M$ is any $\nu$-module, let $\tau$ operates on $M^p$ by letting $\alpha$ operates by cyclic permutation and by letting $\alpha_i$ operate on the $i$-th factor as does $\beta.$ Let $\alpha_i$ operates trivially on $W_1.$ Then $\tau$ operates on $W_1$ and we let $\tau$ operate diagonally on $W_1\otimes M^p.$ In particular $W_1\otimes W_2^p$ is then a $\Lambda\tau$-free resolution of $\Lambda$.

Let $(K,\theta)\in \mathfrak{C}(p).$ We let $\Sigma_{p^2}$ operate on $K^{p^2}$ by permutations with the $(i,j)$-th factor $K$ being the $j$-th factor $K$ in the $i$-th factor $K^p$ of $K^{p^2}=(K^p)^p.$ We let $\nu$ operate on $W_2\otimes K^p$ by letting $\beta$ acting as cyclic permutation on $K^p$. By the previous paragraph this fixes an action of $\tau$ on $W_1\otimes (W_2\otimes K^p)^p.$

Let $Y$ be any $\Lambda\Sigma_{p^2}$-free resolution of $\Lambda$ with $Y_0=\Lambda\Sigma_{p^2}$ and let $w\colon W_1\otimes W_2^p\rightarrow Y$ be any $\tau$-morphism over $\Lambda.$

With these notations, we have the following definition.
\begin{Def}\lb{adem}
 Let $(K,\theta)\in \mathfrak{C}(p).$ We say that $(K,\theta)$ is an Adem object if there exists a $\Sigma_{p^2}$-morphism $\eta\colon Y\otimes K^{p^2}\rightarrow K$ such that the following diagram is $\Lambda\tau$-homotopy commutative\
$$
\begin{diagram}
 \node{(W_1\otimes W_2^p)\otimes K^{p^2}}\arrow{e,t}{w\otimes id}\node{Y\otimes K^{p^2}}\arrow{e,t}{\eta}\node{K}\\
\node{W_1\otimes(W_2\otimes K^p)^p}\arrow{n,l}{id\times \tilde{U}}\arrow{e,t}{id\otimes \theta^p}\node{W_1\otimes K^p}\arrow{e,t}{\theta}\node{K}\arrow{n,r}{id}
\end{diagram}
$$
Here $\tilde{U}$ is the shuffle map and $\Sigma_{p^2}$ acts trivially on $K.$
\end{Def}

It has been proved in \cite{may1} that the following relations among the $\mathcal{P}^s$ and $\beta\mathcal{P}^s$ are valid on all cohomology classes of Adem objects in $\mathfrak{C}(p),~p>2$.
\begin{itemize}
 \item  If $a<pb$ then $\beta^{e}\mathcal{P}^a\mathcal{P}^b=\sum_{i}(-1)^{a+i}(a-pi,(p-1)b-a+i-1)\beta^e\mathcal{P}^{a+b-i}\mathcal{P}^i$.
\item  If $a\leq pb$ then $\beta^e\mathcal{P}^a\beta\mathcal{P}^b=(1-e)\sum_i(-1)^{a+i}(a-pi,(p-1)b-a+i-1)\beta\mathcal{P}^{a+b-i}\mathcal{P}^i\\~~~~~~~~~~~~~~~~~~~~~~~~~~~ ~~~~~~~~ ~ - \sum_i(-1)^{a+i}(a-pi-1,(p-1)b-a+i)\beta^e\mathcal{P}^{a+b-i}\beta\mathcal{P}^i$.
\end{itemize}
where $e=0,1$ and $\beta^0\mathcal{P}^s=\mathcal{P}^s$ and $\beta^1\mathcal{P}^s=\beta\mathcal{P}^s.$

 We apply the above algebraic construction to define Steenrod reduced power operations in equivariant simplicial cohomology of an $O_G$-simplicial set, as defined in the Section $4$. This is done by constructing a functor $\Gamma$ from $\mathcal{A}_{\Lambda}$ to $\mathfrak{C}(p).$

Let $(T,M_0,\rho)$ be an object of $\mathcal{A}_{\Lambda}$. We have already noted in Remark \ref{4.1} that the cochain complex $C_{\rho}^*(T;M_0)$, equipped with the cup product, is a homotopy associative differential graded $\Lambda$-algebra. We now construct a morphism of $\Lambda\pi$-complexes $$\theta\colon W\otimes C_{\rho}^*(T;M_0)^p\rightarrow C_{\rho}^*(T;M_0)$$ so that $(C_{\rho}^*(T;M_0),\theta)$ becomes an object of the category $\mathfrak{C}(p).$

For a simplicial set $L$, let $C(L)$ denote the normalized chain complex with coefficients $\Lambda.$ We recall the following lemma from \cite{may1}.
\begin{Lem}\lb{5.2}
Let $\pi$ be a subgroup of $\Sigma_p$ (not necessarily cyclic of order $p$) and $W$ be a $\Lambda\pi$-free resolution of $\Lambda$ such that $W_0=\Lambda\pi$ with generator $e_0.$ For simplicial sets $L_1,\cdots,L_p$, there exists a morphism in $ch_{\Lambda}$ $$\Phi\colon W\otimes C(L_1\times\cdots\times L_p)\rightarrow W\otimes C(L_1)\otimes \cdots \otimes C(L_p)$$ which is natural in the $L_i$ and satisfies the following properties
\begin{enumerate}
\item For $\sigma\in\pi,$ the following diagram is commutative
$$
 \begin{diagram}
\node{W\otimes C(L_1\times\cdots\times L_p)}\arrow{e,t}{\Phi}\arrow{s,l}{\sigma}\node{W\otimes C(L_1)\otimes \cdots \otimes C(L_p)}\arrow{s,r}{\sigma}\\
\node{W\otimes C(L_{\sigma(1)}\times\cdots\times L_{\sigma(p)})}\arrow{e,t}{\Phi}\node{W\otimes C(L_{\sigma(1)})\otimes \cdots \otimes C(L_{\sigma(p)})}
 \end{diagram}
$$

\item  $\Phi$ is the identity homomorphism on $W\otimes C_0(L_1\times\cdots \times L_p).$
\item $\Phi(e_0\otimes (x_1,\cdots,x_p))=e_0\otimes \xi(x_1,\cdots,x_p)$ where $x_i\in L_j$ for $1\leq i\leq p$ and $$\xi\colon C(L_1\times\cdots\times L_p)\rightarrow C(L_1)\otimes \cdots \otimes C(L_p)$$ is the Alexander-Whitney map.
\item $\Phi(W\otimes C_j(L_1\times\cdots \times L_p))\subseteq\sum_{k\leq pj} W\otimes [C(L_1)\otimes \cdots \otimes C(L_p)]_k$
\item Any two such $\Phi$ are naturally equivariantly homotopic.
\end{enumerate}
\end{Lem}
In the special case $L_1=\cdots =L_p=L,$ we obtain a natural morphism of chain complexes of $\Lambda\pi$-modules $$\Phi\colon W\otimes C(L^p)\rightarrow W\otimes C(L)^p$$ which satisfies the last four conditions of the Lemma \ref{5.2}.

Let $T\in O_G\mathcal{S}$. Applying the above special case of the Lemma \ref{5.2} to each simplicial set $T(G/H)$, we obtain chain maps $\Phi_H\colon W\otimes C(T(G/H)^p)\rightarrow W\otimes C(T(G/H))^p$ which is $\pi$-equivariant. Since $\Phi_H$ is natural with respect to maps of simplicial sets, we see that $\Phi_H\circ (id_W\otimes C(T(\hat{a})^p))=(id_W\otimes C(T(\hat{a}))^p)\circ \Phi_K$ where $a^{-1}Ha\subseteq K$. Thus we have a morphism $\underline{\Phi}$ of $O_G$-chain complexes $$\underline{\Phi}\colon W\otimes \underline{C}(T^p)\rightarrow W\otimes \underline{C}(T)^p,~\mbox{defined by}~\underline{\Phi}(G/H)=\Phi_H,~H\subseteq G.$$

Now suppose that on $O_G$-group $\rho$ operates on $T$. The diagonal action of $\rho$ on $T^p$ induces $\rho$-action on $\underline{C}(T^p)$. Also we have an induced $\rho$-action on $\underline{C}(T)$. We let $\rho$ operate diagonally on $\underline{C}(T)^p$ and trivially on $W$. The naturality of $\Phi_H$ with respect to maps from $T(G/H)$ into itself shows that $\Phi_H$ is $\rho(G/H)$-equivariant. Thus $\underline{\Phi}$ is $(\pi\times \rho)$-equivariant.
Hence we obtain the following.

\begin{Cor} \lb{diag}
Let $T\in O_G\mathcal{S}$ and an $O_G$-group $\rho$ operates on $T$. For a subgroup $\pi$ of $\Sigma_p$ ($\pi$ not necessarily cyclic of order $p$), let $W$ be a $\Lambda\pi$-free resolution of $\Lambda$ such that $W_0=\Lambda\pi$ with generator $e_0.$ Then there is a natural transformation $$\underline{\Phi}\colon W\otimes \underline{C}(T^p)\rightarrow W\otimes \underline{C}(T)^p$$ such that
\begin{enumerate}
\item $\underline{\Phi}$ is $(\pi\times \rho)$-equivariant.
\item  $\underline{\Phi}$ is the identity homomorphism on $W\otimes \underline{C}_0(T^p).$
\item $\underline{\Phi}(G/H)(e_0\otimes (x_1,\cdots,x_p))=e_0\otimes \underline{\xi}(G/H)(x_1,\cdots,x_p)$ where $x_i\in T(G/H)$ for $1\leq i\leq p$ and $\underline{\xi}(G/H)\colon C(T(G/H)^p)\rightarrow C(T(G/H))^p$ is the Alexander-Whitney map.
\item $\underline{\Phi}(G/H)(W\otimes C_j(T(G/H)^p))\subseteq\sum_{k\leq pj}W\otimes (C(T(G/H))^p)_k$
\item The map $\underline{\Phi}$ is natural with respect to equivariant maps of $O_G$-simplicial sets and any two such $\underline{\Phi}$ are naturally equivariantly homotopic.
\end{enumerate}
\end{Cor}

\begin{Def}\lb{definition}
For $(T,M_0,\rho)\in \mathcal{A}_{\Lambda}$, let $D\colon T\rightarrow T^p$ be the diagonal map $$D(G/H)(x)=(x,\cdots ,x), ~x\in T(G/H),$$ inducing a map $D_*\colon \underline{C}(T)\rightarrow \underline{C}(T^p).$ Let $\underline{\Delta}\colon W\otimes \underline{C}(T)\rightarrow \underline{C}(T)^p$ to be the composite
$$\underline{\Delta}\colon W\otimes \underline{C}(T)\xrightarrow{id\otimes D_*} W\otimes\underline{C}(T^p)\xrightarrow{\underline{\Phi}}W\otimes \underline{C}(T)^p\rightarrow \underline{C}(T)^p$$ where the last map is the augmentation. Observe that the map $\underline{\Delta}$ is $(\pi\times \rho)$-equivariant. Moreover we have a natural map $$\alpha\colon  [C_{\rho}^*(T;M_0)]^p\rightarrow Hom_{\rho}(\underline{C}(T)^p,M_0)$$ defined by $$\alpha(f_1\otimes \cdots \otimes f_p)(G/H)(x_1\otimes \cdots \otimes x_p)=(-1)^a f_1(G/H)(x_1)\cdots f_p(G/H)(x_p),$$ where $f_i\in C_{\rho}^*(T), x_i\in \underline{C}(T)(G/H), i=1,\cdots,p$ and $a=\Pi_{k=1}^p deg(x_k).$ Hence dualising $\underline{\Delta}$ we get a natural morphism of $\Lambda\pi$-complexes
$$\theta\colon W\otimes C^*_{\rho}(T;M_0)^p\rightarrow C^*_{\rho}(T;M_0)$$ $$\theta(w\otimes f)(G/H)(x)=(-1)^{deg(w)deg(x)}\alpha(f)(G/H)(\underline{\Delta}(G/H)(w\otimes x))$$ where $w\in W, f\in C^*_{\rho}(T)^p,x\in C(T(G/H)).$

Note that $\theta(e_0\otimes f)= D^*\xi^*\alpha(f)$ for any $f\in C_{\rho}^*(T)^p.$ As before let $V$ denote a $\Lambda\Sigma_p$-free resolution of $\Lambda$ and $j\colon W\rightarrow V$ be the map induced by the inclusion $\pi\hookrightarrow \Sigma_p.$ We apply the Corollary \ref{diag} for the (sub)group $\Sigma_p$ to get $\tilde{\Phi}\colon V\otimes \underline{C}(T^p)\rightarrow W\otimes \underline{C}(T)^p$. Then $\tilde{\Phi}\circ(j\otimes id)$ satisfies first four conditions of the Corollary \ref{diag} for the subgroup $\pi$ and hence must be equivariantly homotopic to $\Phi.$ Therefore $\tilde{\theta}\colon V\otimes C^*_{\rho}(T)^p\rightarrow C^*_{\rho}(T)$ can be defined such that $\tilde{\theta}\circ (j\otimes id)$ is $\Lambda \pi$-equivariantly homotopic to $\theta.$  Therefore $(C^*_{\rho}(T),\theta)$ is an object of the category $\mathfrak{C}(p).$ Thus we obtain a contravariant functor $\Gamma\colon \mathcal{A}_{\Lambda}\rightarrow \mathfrak{C}(p)$ by letting $\Gamma(T,M_0,\rho)=(C^*_{\rho}(T;M_0),\theta)$ and $\Gamma(f_0,f_1,f_2)= C^*(f_0,f_1,f_2)$ on morphisms (see Remark \ref{map}).
\end{Def}

The next lemma is the key to show that $(C^*_{\rho}(T;M_0),\theta)$ is a Cartan object of $\mathfrak{C}(p).$
Let $\phi=(\epsilon \otimes id)\Phi$ where $\Phi$ is obtained from the Lemma \ref{5.2} and $\epsilon\colon W\rightarrow \Lambda$ is the augmentation.

\begin{Lem}\lb{lem 5.5}
Let $L_i,S_i~i=1,\cdots, p$ be simplicial sets. Let $u\colon (\prod_{i=1}^{p}L_i\times \prod_{i=1}^{p}S_i)\rightarrow \prod_{i=1}^{p} (L_i\times S_i)$ and $U\colon (\otimes_{i=1}^p C(L_i))\otimes (\otimes_{i=1}^p C(S_i))\rightarrow \otimes_{i=1}^{p}[C(L_i)\otimes C(S_i)]$ be shuffle maps. Let $t$ denote the flip map, that is, $t(x\otimes y)=y\otimes x$. Then there exists a homotopy $$\mathcal{H}\colon W\otimes C(\prod_{i=1}^{p}L_i\times \prod_{i=1}^{p}S_i)\rightarrow \otimes_{i=1}^{p}[C(L_i)\otimes C(S_i)]$$ of the chain maps $\xi^p \phi(id\otimes u)$ and ${U(\phi\otimes \phi) (id\otimes t\otimes id)(\psi\otimes id \otimes id)(id\times \xi)}$, so that the following diagram is homotopy commutative.
$$
\begin{diagram}
 \node{W\otimes C(\prod_{i=1}^{p}L_i\times \prod_{i=1}^{p}S_i)}\arrow{e,t}{id\times u}\arrow{s,l}{id\otimes \xi}\node{W\otimes C(\Pi_{i=1}^{p} (L_i\times S_i))}\arrow{e,t}{\phi}\node{\otimes_{i=1}^{p}[C(L_i\times S_i)]}\arrow{s,r}{\xi^p}\\
\node{W\otimes C(\prod_{i=1}^{p}L_i)\otimes C(\prod_{i=1}^{p}S_i)}\arrow[2]{e,b}{U(\phi\otimes \phi)(id\otimes t\otimes id)(\psi\otimes id \otimes id)}\node[2]{\otimes_{i=1}^{p}[C(L_i)\otimes C(S_i)]}
\end{diagram}
$$

Moreover the homotopy $\mathcal{H}$ is natural in the $L_i,S_i$ and following diagram commutes for $\sigma\in \pi.$
$$
\begin{diagram}
\node{W\otimes C(\prod_{i=1}^{p}L_i\times \prod_{i=1}^{p} S_i)}\arrow{e,t}{\mathcal{H}}\arrow{s,l}{\sigma\otimes \sigma}\node{\otimes_{i=1}^{p}[C(L_i)\otimes C(S_i)]}\arrow{s,r}{\sigma}\\
\node{W\otimes C(\prod_{i=1}^{p}L_{\sigma(i)}\times \prod_{i=1}^{p}S_{\sigma(i)})}\arrow{e,b}{\mathcal{H}}\node{\otimes_{i=1}^{p}[C(L_{\sigma(i)})\otimes C(S_{\sigma(i)})]}
\end{diagram}
$$
\end{Lem}

\begin{proof}
The proof is similar to the proof of Lemma $7.1$ of \cite{may1}.
Let $A_j=C_j(\prod_{i=1}^{p}L_i\times \prod_{i=1}^{p}S_i)$ and $B_j=[\otimes_{i=1}^{p}C(L_i)\otimes C(S_i)]_j.$ We construct $\mathcal{H}$ on $W_i\otimes A_j$ by induction on $i$ and for fixed $i$ by induction on $j.$ Note that the two maps agree on $W\otimes A_0.$ So $H$ is the zero map on $W\otimes A_0.$ To define $\mathcal{H}$ on $W_0\otimes A_j,~j\geq 0,$ it suffices to define on $e_0\otimes A_j$, since $\mathcal{H}$ can then be uniquely extended to all of $W_0\otimes A_j$ using the commutativity of the second diagram. The functor $e_0\otimes A_j$ is represented by the model $\Delta[j]^p\times \Delta[j]^p$ and $W\otimes B_j$ is acyclic on this model. Therefore, by acyclic model argument, $\mathcal{H}$ can be defined on $e_0\otimes A_j$, provided $\mathcal{H}$ is known on $e_0\otimes A_{j-1}.$ But $\mathcal{H}$ has already been defined on $W_0\otimes A_0.$ Hence by induction on $j$, we can define $\mathcal{H}$ on $e_0\otimes A_j,~j\geq 0.$ To define $\mathcal{H}$ on $W_i\otimes A_j,$ assume that it has already been defined on $W_{i^{\prime}}\otimes A_j,~i^{\prime}<i,~j\geq 0$ and on $W_i\otimes A_{j^{\prime}},~j^{\prime}<j.$ Choose a $\Lambda\pi$-basis $\{w_k\}$ for $W_i$. As before, it suffices to define $\mathcal{H}$ on $w\otimes A_j,~w\in\{w_k\}.$ We can repeat the acyclic model argument replacing $e_0$ by $w$, and hence we are through by induction.
\end{proof}

In the special case $L_1=\cdots=L_p=L,~S_1=\cdots =S_p=S,$ we obtain the following corollary.
\begin{Cor}\lb{cor 5.6}
 For simplicial sets $L,S$ the two chain maps $\xi^p \phi(id\otimes u)$ and \linebreak ${U(\phi\otimes \phi) (id\otimes t\otimes id)(\psi\otimes id \otimes id)(id\times \xi)}$ from
$W\otimes C(L^p\times S^p)$ to $[C(L)\otimes C(S)]^p$
are $\Lambda\pi$-equivariantly homotopic and the homotopy is natural in $L$ and $S$.
\end{Cor}

Suppose $(T,M_0,\rho)$ and $(T^{\prime},M_0^{\prime},\rho^{\prime})$ are objects of $\mathcal{A}_{\Lambda}.$ For the product actions of $\rho\times \rho^{\prime}$ on $T\times T^{\prime}$ and $M_0\otimes M_0^{\prime}$, we have $(T\times T^{\prime},M_0 \times M_0^{\prime},\rho\times \rho^{\prime})\in \mathcal{A}_{\Lambda}$. The following lemma relates $\Gamma(T\times T^{\prime},M_0 \times M_0^{\prime},\rho\times \rho^{\prime})=(C_{\rho\times \rho^{\prime}}^*(T\times T^{\prime};M_0\otimes M_0^{\prime}),\theta)$ to $$\Gamma(T,M_0,\rho)\otimes \Gamma(T^{\prime},M_0^{\prime},\rho^{\prime})=(C_{\rho}^*(T;M_0)\otimes C_{\rho^{\prime}}^*(T^{\prime};M_0^{\prime}),\tilde{\theta}).$$ Let $$\tilde{\alpha}\colon C_{\rho}^*(T;M_0)\otimes C_{\rho^{\prime}}^*(T^{\prime};M_0^{\prime})\rightarrow Hom_{\rho\times \rho^{\prime}}(\underline{C}(T)\otimes \underline{C}(T^{\prime}),M_0\otimes M_0^{\prime})$$ be defined as
$$\tilde{\alpha}(f\otimes g)(G/H)(x\otimes y)=(-1)^{deg(x)deg(y)}f(G/H)(x)\otimes g(G/H)(y), ~ H\subseteq G,$$
where $f\in C_{\rho}^*(T;M_0),~g\in C_{\rho^{\prime}}^*(T^{\prime};M_0^{\prime}),~x\in \underline{C}(T)(G/H),~y\in \underline{C}(T^{\prime})(G/H).$
\begin{Lem}\lb{lem 5.7}
With notations as above, the following diagram is $\Lambda\pi$-homotopy commutative.
$$
\begin{diagram}
\node{W\otimes C_{\rho\times \rho^{\prime}}^*(T\times T^{\prime};M_0\otimes M_0^{\prime})^p}\arrow{e,t}{\theta}\node{C_{\rho\times \rho^{\prime}}^*(T\times T^{\prime};M_0\otimes M_0^{\prime})} \\
\node{W\otimes [C_{\rho}^*(T;M_0)\otimes C_{\rho^{\prime}}^*(T^{\prime};M_0^{\prime})]^p}\arrow{n,l}{id\otimes(\underline{\xi}^*\tilde{\alpha})^p}\arrow{e,t}{\tilde{\theta}}\node{C_{\rho}^*(T;M_0)\otimes C_{\rho^{\prime}}^*(T^{\prime};M_0^{\prime})}\arrow{n,r}{\underline{\xi}^*\tilde{\alpha}}
\end{diagram}
$$
\end{Lem}
\begin{proof}
Let $D,~D^{\prime},~\tilde{D}$ be the diagonals for $T,~T^{\prime},~T\times T^{\prime}$ respectively. Let $\underline{u}\colon T^p\times T^{\prime p}\rightarrow (T\times T^{\prime})^p$ and $\underline{U}\colon \underline{C}(T)^p\otimes \underline{C}(T^{\prime})^p\rightarrow [\underline{C}(T)\otimes \underline{C}(T^{\prime})]^p$ be the shuffle maps. Let $t$ be the switch map.

 By definition of $\theta$ and $\tilde{\theta},$ it suffices to prove that the following diagram of $O_G$-chain complexes is $\Lambda(\pi\times \rho\times \rho^{\prime})$-equivariant homotopy commutative.

$${(A)\cdots}
\begin{diagram}
\node{W\otimes {\underline{C}(T\times T^{\prime})}}\arrow{e,t}{\underline{\Delta}}\arrow{s,l}{id\times \underline{\xi}}\node{\underline{C}(T\times T^{\prime})^p}\arrow{s,r}{\underline{\xi}^p}\\
\node{W\otimes \underline{C}(T)\otimes \underline{C}(T^{\prime})}\arrow{e,t}{\zeta}\node{[\underline{C}(T)\otimes \underline{C}(T^{\prime})]^p}
\end{diagram}
$$

Here $$\underline{\Delta}=(\epsilon\otimes id)\underline{\Phi}(id\otimes \tilde{D}) ,~\zeta =\underline{U}(\underline{\Delta}\otimes \underline{\Delta})(id\otimes t\otimes id)(\psi\otimes id\otimes id).$$ Let $\underline{\phi}=(\epsilon\otimes id)\underline{\Phi}.$ Observe that $\tilde{D}=\underline{u}(D\times D^{\prime})$ and $$(id\otimes D\otimes id \otimes D^{\prime})(id\otimes t\otimes id)(\psi\otimes id\otimes id)=(id\otimes t\otimes id)(\psi\otimes id\otimes id)(id\otimes D\otimes D^{\prime}).$$
Observe that the following diagram commutes by the naturality of $\underline{\xi}$.

$${(A1)\cdots}
\begin{diagram}
\node{W\otimes \underline{C}(T\times T^{\prime})}\arrow{e,t}{id\otimes (D\times D^{\prime})}\arrow{s,l}{id\otimes \underline{\xi}}\node{W\otimes \underline{C}(T^p\times T^{\prime p})}\arrow{s,l}{id\otimes \underline{\xi}}\\
\node{W\otimes \underline{C}(T)\otimes \underline{C}(T^{\prime})}\arrow{e,t}{id\otimes D \otimes D^{\prime}}\node{W\otimes \underline{C}(T^p)\otimes \underline{C}(T^{\prime p})}
\end{diagram}
$$
Let $\mathcal{F}$ denote the following diagram of $O_G$-chain complexes of $\Lambda$-modules.
$${(A2)\cdots}
\begin{diagram}
 \node{W\otimes \underline{C}(T^p\times T^{\prime p})}\arrow{e,t}{id\otimes \underline{u}}\arrow{s,l}{id\otimes \underline{\xi}}\node{W\otimes \underline{C}([T\times T^{\prime}]^p)}\arrow{e,t}{\underline{\phi}}\node{\underline{C}(T\times T^{\prime})^p}\arrow{s,r}{\underline{\xi}^p}\\
\node{W\otimes \underline{C}(T^p)\otimes \underline{C}(T^{\prime p})}\arrow[2]{e,t}{\underline{U}(\underline{\phi}\otimes \underline{\phi})(id\otimes t\otimes id)(\psi\otimes id\otimes id)}\node[2]{[\underline{C}(T)\otimes \underline{C}(T^{\prime})]^p}
\end{diagram}
$$
  Then $\mathcal{F}(G/H)$ is $\Lambda\pi$-homotopy commutative, by Corollary \ref{cor 5.6}. The naturality of this homotopy with respect to maps of $T(G/H)$ into itself implies that the homotopy is equivariant for the $\rho(G/H)$-action on $T(G/H).$ Similarly it is $\rho^{\prime}(G/H)$-equivariant. These natural equivariant homotopies of chain complexes combine together to form $\Lambda(\pi\times \rho\times \rho^{\prime})$-equivariant homotopy, which makes the diagram $(A2)$ $\Lambda(\pi\times \rho\times \rho^{\prime})$-equivariant homotopy commutative.

Now observe that the diagram $(A)$ is juxtaposition of the diagrams $(A1)$ and $(A2)$. Hence the diagram $(A)$ is $\Lambda(\pi\times \rho\times \rho^{\prime})$-equivariant homotopy commutative.
\end{proof}

\begin{Prop}
 For an object $(T,M_0,\rho)$ of $\mathcal{A}_{\Lambda}$, $\Gamma(T,M_0,\rho)=(C_{\rho}^*(T;M_0),\theta)$ is a Cartan object of $\mathfrak{C}(p).$
\end{Prop}
\begin{proof}
Recall that $(C_{\rho}^*(T;M_0),\theta)$ is called a Cartan object if the cup product is a morphism of $\mathfrak{C}(p).$
Now observe that $(T,M_0,\rho)\xrightarrow{(D,id,id)} (T\times T,M_0,\rho)\xrightarrow{(id,m,D)} (T\times T,M_0\otimes M_0,\rho\times \rho)$ are morphisms in $\mathcal{A}_{\Lambda},$ where $m\colon M_0\otimes M_0\rightarrow M_0$ is the multiplication, $D$ denotes the diagonal map, and we let $\rho$ operates diagonally on $T\times T$.

 Applying the Lemma \ref{lem 5.7} with $(T,M_0,\rho)=(T^{\prime},M_0^{\prime},\rho^{\prime}),$ and composing with the morphism $C^*(id,m,D)$, we see that the composite $\underline{\xi}^*\alpha$ $$C_{\rho}^*(T;M_0)\otimes C_{\rho}^*(T;M_0)\xrightarrow{\alpha} Hom_{\rho}(\underline{C}(T)\otimes \underline{C}(T), M_0)\xrightarrow{\underline{\xi}^*} C_{\rho}^*(T\times T;M_0)$$
is a morphism in $\mathfrak{C}(p).$ Also $C^*(D,id,id)\colon C_{\rho}^*(T\times T;M_0)\rightarrow C_{\rho}^*(T;M_0)$ is a morphism in $\mathfrak{C}(\Lambda).$ Hence the cup product is a morphism in $\mathfrak{C}(p).$
\end{proof}

Next we show that $(C_{\rho}^*(T;M_0)$ is an `Adem object' in $\mathfrak{C}(p).$
\begin{Prop}
 For $(T,M_0,\rho)\in \mathcal{A}_{\Lambda},$ $\Gamma(T,M_0,\rho)=(C_{\rho}^*(T;M_0),\theta)$ is an Adem object of $\mathfrak{C}(p).$
\end{Prop}
\begin{proof}
 With the notations of the Definition \ref{adem}, we first construct $$\eta\colon Y\otimes C_{\rho}^*(T;M_0)^{p^2}\rightarrow C_{\rho}^*(T;M_0).$$ The procedure is similar to the construction of $\theta$. We remark that the proof of the Lemma \ref{5.2} works for any subgroup $\pi$ of $\Sigma_r$, $r$ being any positive integer. Thus, for simplicial sets $L_1,\cdots, L_r$, we have a chain map $$\Phi\colon Y\otimes C(L_1\times \cdots L_r)\rightarrow Y\otimes C(L_1)\otimes\cdots \otimes C(L_r),$$ satisfying properties of the Lemma \ref{5.2}. Specializing to $L_1=\cdots =L_r=L$ and $\pi=\Sigma_r$, and then passing to $O_G$-simplicial set $T$ with $O_G$-group action $\rho$, we obtain $O_G$-chain map $\underline{\Delta}\colon Y\otimes \underline{C}(T)\rightarrow \underline{C}(T)^{p^2}$ which is $(\Sigma_r\times \rho)$-equivariant. Proceeding similarly to the construction of the map $\theta$, we obtain $\eta$.

 Dualising the diagram in Definition \ref{adem}, it suffices to prove that the following diagram is $\Lambda(\tau\times \rho)$-homotopy commutative,
$$
\begin{diagram}
 \node{W_1\otimes W_2^p\otimes \underline{C}(T)}\arrow{e,t}{w\otimes id}\arrow{s,l}{t\times id}\node{Y\otimes \underline{C}(T)}\arrow{e,t}{\underline{\Delta}}\node{\underline{C}(T)^{p^2}}\\
\node{W_2^p\otimes W_1\otimes \underline{C}(T)}\arrow{e,t}{id\otimes \underline{\Delta}}\node{W_2^p\otimes \underline{C}(T)^p}\arrow{e,t}{\underline{U}}\node{[W_2\otimes \underline{C}(T)]^p}\arrow{n,r}{\underline{\Delta}^p}
\end{diagram}
$$
where the notations are as in Lemma \ref{lem 5.7}. Then, as in \cite{may1}, it suffices to show that the the natural transformations $\chi,\Omega\colon W_1\otimes W_2^p\otimes \underline{C}(T^{p^2})\rightarrow \underline{C}(T)^{p^2}$, defined by $\chi=\underline{\phi}(w\otimes id)$ and $\Omega=\underline{\phi}^p\underline{U}(id\otimes \underline{\phi})(t\otimes id)$, are $\Lambda(\tau\times \rho)$-equivariantly homotopic. Here $\tau$ operates by permutation of factors and the $O_G$-group $\rho$ operates diagonally on $T^{p^2}$ and $\underline{C}(T)^{p^2}.$ Now for simplicial sets without any group action, replacing $T^{p^2}$ by the product $\prod_{i,j=1}^p L_{ij}$ and then following the methods of the Lemma \ref{lem 5.5}, the corresponding chain maps can be shown to be $\tau$-equivariantly homotopic and the homotopy is natural with respect to maps of simplicial sets. Hence $\chi$ and $\Omega$ are $\Lambda(\tau\times \rho)$-equivariantly homotopic.
\end{proof}

Thus we have the following theorem.
\begin{Thm}\lb{Th 5.12}
 Let $(T,M_0,\rho)\in \mathcal{A}_{\Lambda}$, $\Lambda=\mathbb{Z}_p,~p>2.$ Then there exist functions $$\mathcal{P}^s\colon H^q_{\rho}(T;M_0)\rightarrow H^{q+2s(p-1)}_{\rho}(T;M_0),$$ $$\beta\mathcal{P}^s\colon H^q_{\rho}(T;M_0)\rightarrow H^{q+2s(p-1)+1}_{\rho}(T;M_0),$$ which satisfies the following properties
\begin{enumerate}
\item $\mathcal{P}^s$ and $\beta\mathcal{P}^s$ are natural homomorphisms.
 \item $\mathcal{P}^s=\beta \mathcal{P}^s=0$ if $s<0.$ Also $\mathcal{P}^s(x)=0$ if $2s>q$, $\beta\mathcal{P}^s=0$ if $2s\geq q$.
\item  $\mathcal{P}^s(x)=x^p$ if $2s=q.$
\item {(Cartan formula)} For $x,y\in H^q_{\rho}(T;M_0)$, $$\mathcal{P}^s(x\cup y)=\sum_{i+j=s}\mathcal{P}^i(x)\cup \mathcal{P}^j(y),$$
$$\beta\mathcal{P}^{s+1}(x\cup y)=\sum_{i+j=s}[\beta\mathcal{P}^{i+1}(x)\cup \mathcal{P}^j(y)+(-1)^{deg(x)}\mathcal{P}^i(x)\cup \beta\mathcal{P}^{j+1}(y)].$$
\item {(Adem relation)} If $a<pb$ then $$\beta^{e}\mathcal{P}^a\mathcal{P}^b=\sum_i(-1)^{a+i}(a-pi,(p-1)b-a+i-1)\beta^e\mathcal{P}^{a+b-i}\mathcal{P}^i.$$
If $a\leq pb$ then $$\beta^e\mathcal{P}^a\beta\mathcal{P}^b=(1-e)\sum_i(-1)^{a+i}(a-pi,(p-1)b-a+i-1)\beta\mathcal{P}^{a+b-i}\mathcal{P}^i$$ $$~~~~~~~~~~~~~~~~~~~-\sum_i(-1)^{a+i}(a-pi-1,(p-1)b-a+i)\beta^e\mathcal{P}^{a+b-i}\beta\mathcal{P}^i,$$
where $e=0,1$ and $\beta^0\mathcal{P}^s=\mathcal{P}^s$ and $\beta^1\mathcal{P}^s=\beta \mathcal{P}^s.$
\end{enumerate}
\end{Thm}
\begin{proof}
 We only need to prove that $\mathcal{P}^s=\beta\mathcal{P}^s=0$ for $s<0$. By the definition of the power operations, it suffices to show that $D_i(x)=0$ for $i>pq-q$, $deg(x)=q$. Recall that $\underline{\Delta}=(\epsilon\otimes id)\underline{\Phi}(id\times D)$ and $$\underline{\Phi}(e_i\otimes D(k))\in \sum_{j<pq}W_{pq-j}\otimes [\underline{C}(T)]_j^p\subseteq Ker(\epsilon\otimes id)~\mbox{for}~i>pq-q.$$ Hence $\underline{\Delta}(e_i\otimes k)=0$ for $k\in \underline{C}_{pq-i}(T)$.
\end{proof}

Let $X$ be a one vertex $G$-Kan complex and $M$ be an equivariant local coefficients on $X$. We define the reduced power operations in the simplicial Bredon-Illman cohomology with local coefficients by $\mathcal{P}^s=\mu^{* -1}\mathcal{P}^s\mu^*$ and $\beta\mathcal{P}^s=\mu^{*-1}(\beta\mathcal{P}^s)\mu^*$, where the symbols $\mathcal{P}^s$ and $\beta\mathcal{P}^s$ on the right side of the equality denote the power operations as constructed in the category $\mathcal{A}_{\Lambda}$ and $\mu^*$ is the isomorphism as given in the Theorem \ref{eilenberg}. Thus we have the following theorem.

\begin{Thm} Let $X$ be a one vertex $G$-Kan complex and $M$ be an equivariant local coefficients on $X$.
 Then there exist natural homomorphisms $$\mathcal{P}^s\colon H^q_G(X;M)\rightarrow H^{q+2s(p-1)}_G(X;M),$$ $$\beta\mathcal{P}^s\colon H^q_G(X;M)\rightarrow H^{q+2s(p-1)+1}_G(X;M),$$ which satisfies the properties $(1)-(5)$ of the Theorem \ref{Th 5.12}.

If $G$ is trivial then $\mathcal{P}^s$ can be naturally identified with the reduced power operations in local coefficients \cite{git}.
\end{Thm}

\begin{proof}
Since the isomorphism $\mu^*$ of Eilenberg theorem is natural and respects the cup product, the first part follows from the Theorem \ref{Th 5.12}.

 For the second part, we just remark that when $G$ is trivial the map $\underline{\Delta}$, as constructed in the Definition \ref{definition}, reduces to the $(\pi\times \rho)$-equivariant chain mapping $\phi^{\prime}$ of \cite{git} (see Section $4.2$ of \cite{git}).
\end{proof}


{\bf Goutam Mukherjee}\\
Indian Statistical Institute, Kolkata-700108, India.\\
e-mail: goutam@isical.ac.in

{\bf Debasis Sen}\\
Indian Statistical Institute, Kolkata-700108, India.\\
e-mail: dsen\_r@isical.ac.in

\end{document}